\documentclass[11pt,titlepage]{article}
\usepackage[margin=1.2in]{geometry}
\usepackage{amsmath, amsbsy, amsthm, array, mathrsfs}
\usepackage{physics}
\usepackage{amssymb,latexsym,verbatim}
\usepackage{color}
\usepackage{dsfont, bbm}
\usepackage{relsize}
\usepackage{subcaption}

\newcommand{\R}{\mathbb{R}}

\newcommand{\noin}{\noindent}
\newcommand{\bee}{\begin{eqnarray*}}
\newcommand{\ene}{\end{eqnarray*}}
\newcommand{\bec}{\begin{center}}
\newcommand{\enc}{\end{center}}
\newcommand{\be}{\begin{equation}}
\newcommand{\ee}{\end{equation}}

\newcommand{\ep}{\varepsilon}
\newcommand{\mb}{\mathbf}
\newcommand{\bs}{\boldsymbol}
\newcommand{\tb}{\textbf}
\newcommand{\pend}{$\square$}
\newcommand{\vs}{\vskip 3mm}
\newcommand{\bi}{\begin{itemize}}
\newcommand{\ei}{\end{itemize}}

\begin{document}


\title{\LARGE 
Large sample behavior of the least trimmed squares estimator} 
\vs
\vs

\author{ {\sc 
Yijun Zuo}\\[3ex]
         {\small Department of Statistics and Probability} \\[2ex]
         {\small Michigan State University, East Lansing, MI 48824} \\[2ex]
         {\small 
         zuo@msu.edu}\\[6ex]
     }
 \date{\today}
\maketitle

\vskip 3mm
{\small

\begin{abstract}
The least trimmed squares (LTS)  estimator is popular in location, regression, machine learning, and AI literature.	
Despite the empirical version of least trimmed squares (LTS) 
being repeatedly studied in the literature, the population version of the LTS has never been introduced and studied. The lack of the population version hinders the study of the large sample properties of the LTS utilizing the empirical process theory.\vs
\noin
Novel properties of the objective function in both empirical and population settings of the LTS and other properties, are established for the first time in this article.
The primary properties of the objective function
facilitate the establishment of other original results, including the influence function and Fisher consistency. The strong consistency is established with the help of a generalized Glivenko-Cantelli Theorem over a class of functions for the first time.
Differentiability and stochastic equicontinuity promote the establishment of asymptotic normality with a concise and novel approach.

\bigskip
\noindent{\bf AMS 2000 Classification:} Primary 62J05, 62G36; Secondary
62J99, 62G99
\bigskip
\par

\noindent{\bf Key words and phrase:}  trimmed squares of residuals, continuity and differentiability of objective function,  influence function, Fisher consistency, asymptotics. 
\bigskip
\par
\noindent {\bf Running title:} Asymptotics of the least trimmed squares.
\end{abstract}
}
\setcounter{page}{1}

\section {Introduction}

In  classical multiple linear regression analysis, it is assumed that there is a relationship for a given data set $\{(\bs{x}^{\top}_i, y_i)^{\top}, i\in \{1,\cdots, n\}\}$:\vspace*{-2mm}
\be
y_i=(1,\bs{x}^{\top}_i)\bs{\beta}_0+{e}_i,~~ i\in \{1,\cdots, n\},  \label{model.eqn}
\ee
where $y_i$ and $e_i$ (an error term, a random variable, and is assumed to have a zero mean and unknown variance $\sigma^2$ in the classic regression theory) are in $\R^1$, ${\top}$ stands for the transpose, $\bs{\beta}_0=(\beta_{01}, \cdots, \beta_{0p})^{\top}$, the true unknown parameter, {and}~ $\bs{x_i}=(x_{i1},\cdots, x_{i(p-1)})^{\top}$ is in $\R^{p-1}$ ($p\geq 2$) and could be random.  It is seen that $\beta_{01}$ is the intercept term. Writing $\bs{w}_i=(1,\bs{x}^{\top}_i)^{\top}$, one has $y_i=\bs{w}^{\top}_i\bs{\beta}_0+e_i$. The classic assumptions such as linearity and homoscedasticity are implicitly assumed here. Others will be assumed later when they are needed.
\vs
One wants to estimate the $\bs{\beta}_0$ based on the given sample $\mb{z}^{(n)}
:=\{(\bs{x}^{\top}_i, y_i)^{\top}, i\in \{1,\cdots, n\}\}$. (hereafter it is implicitly assumed that they are i.i.d. from parent $(\bs{x}, y)$)
For a candidate coefficient vector $\bs{\beta}$, call the difference between $y_i$ (observed) and $\bs{w^{\top}_i}{\bs{\beta}}$ (predicted),
 the ith residual, $r_i(\bs{\beta})$, ($\bs{\beta}$ is often suppressed). That is,\vspace*{-2mm}
\be r_i:= {r}_i(\bs{\beta})=y_i-\bs{w^{\top}_i}{\bs{\beta}}.\label{residual.eqn}
\ee
To estimate $\bs{\beta}_0$, the classic \emph{least squares} (LS) minimizes  the sum of squares of residuals,
$$\widehat{\bs{\beta}}_{ls}:=\arg\min_{\bs{\beta}\in\R^p} \sum_{i=1}^n r^2_i. $$
Alternatively, one can replace the square above by the absolute value to obtain the least absolute deviations  estimator (aka, $L_1$ estimator, in contrast to the $L_2$ (LS) estimator).
\vs
 Due to its great computability and optimal properties when the error $e_i$ follows a Gaussian 
  distribution, the LS estimator is  popular  in practice across multiple disciplines.
It, however, can behave badly when the error distribution  slightly departs from Gaussian assumption,
particularly when the errors are heavy-tailed or contain outliers.
In fact, both $L_1$ and $L_2$  estimators have the worst $0\%$ asymptotic breakdown point, in sharp contrast to the $50\%$ of the least
trimmed squares estimator (Rousseeuw and Leroy (1987) (RL87)). 
 The latter is one of the most robust alternatives to the LS estimator.
Robust alternatives to the LS estimator are abundant in the literature.
The most popular are, M-estimators (Huber (1964)), 
least median squares (LMS) and least
trimmed squares (LTS) estimators (Rousseeuw (1984) (R84)), 
S-estimators (Rousseeuw and Yohai (1984)), 
MM-estimators (Yohai (1987)), 
$\tau$ -estimators (Yohai and Zamar (1988)), 
 and maximum depth estimators (Rousseeuw and Hubert (1999)), 
 and Zuo (2021 a, b), Zuo and Zuo (2023) (ZZ23), 
 among others.
\vs

 Due to its cube-root consistency of LMS in R84 and its other drawbacks, 
 LTS is preferred over LMS  (see Rousseeuw and Van Driessen (2006) (RVD06)). 
LTS is popular in the literature in view of its fast computibility and high robustness and often severing as the initial estimator for many high breakdown point iterative procedures (e.g., S- and MM- estimators).
The LTS is defined as
the  minimizer of the sum of $h$ \emph{trimmed squares} of residuals. Namely, \vspace*{-2mm}
\be
\widehat{\bs{\beta}}_{lts}:=\arg\min_{\bs{\beta}\in \R^p} \sum_{i=1}^h r^2_{i:n},\label{lts.eqn}
\ee
\vspace*{-.5mm}

\noin
where $r^2_{1:n}\leq r^2_{2:n}\leq \cdots\leq r^2_{n:n}$ are the ordered squared residuals and $\lceil n/2\rceil\leq h < n$ and $\lceil x \rceil$ is the ceiling function.
\vs

There are copious studies on the LTS in the literature.  
Most focused on its computation, e.g., RL87, 
 Stromberg (1993), 
  Hawkins (1994), 
  H\"{o}ssjer (1995), 
  RVD99, RVD06, 
  Hawkins and Oliver (1999), 
  Agull\'{o} (2001), 
  Hofmann et al. (2010), 
 and Klouda (2015). \vs 
Others sporadically addressed the asymptotics,  e.g., RL87, and  Ma\v{s}\'{i}\v{c}ek (2004) 
 addressed the asymptotic normality of the LTS, but limited to the location case, that is, when $p=1$.
 V\'{i}\v{s}{e}k (2006a b c) 
also addressed the asymptotics of the LTS 
 without employing advanced technical tools in a series (three) of lengthy articles for the consistency, root-n consistency, and
asymptotic normality, respectively. The analysis is technically demanding and with  difficult verified 
 assumptions $\mathcal{A, B, C}$. 
Furthermore, 
the analysis is limited to the non-random vectors $\bs{x}_i$s case. In this article, without those assumptions and limitation, those results are established in a concise manner with the help of advanced empirical process theory.
\vs 

 Replacing $(1,\bs{x}^{\top}_i)\bs{\beta}_0$ by a unspecified  nonlinear  function $h(\bs{x}_i, \bs{\beta}_0)$, Chen, at al. (1997) 
 (CSZ97) 
 and C\'{i}\v{z}ek (2004, 2005) (C04, C05) 
  discussed the asymptotics of the LTS in
 a nonlinear regression setting. Now that more general non-linear case has been addressed,  one might wonder is there any merits to discuss the special linear case in this article?\vs
 
 

 There are at least these merits: (i) the nonlinear function $h(\bs{x}_i, \bs{\beta}_0)$ cannot always cover the linear case of $(1,\bs{x}^{\top}_i)\bs{\beta}_0$ for the usual LTS (e.g., in the exponential and power regression cases); (ii) many assumptions for the nonlinear case (see  A1, A2, A3, A4 in CSZ97; H1, H2, H3, H4, H5, H6; D1, D2; I1, I2 in C04 and C05 )  (which are usual difficult to verify) can be dropped for the linear case as demonstrated in this article.
 (iii) A key assumption that  $\{h(\bs{x}, \bs{\beta}), \bs{\beta}\in \Theta\}$ form a VC class of functions over  a compact parameter space $\Theta$ (see CSZ97, C04, C05) can be verified directly in this article.
\vs
To avoid all the drawbacks and limitations above and take advantage of the standard results of the empirical process theory, this article defines the population version of the LTS (Sec.\ref{Sec.2.1}), introduces the novel partition of the parameter space (Sec.\ref{Sec.2.2}), and investigates the
primary properties of the objective function for the LTS both in the empirical and population settings 
(Sec.\ref{Sec.2}) for the first time.  The obtained novel results facilitate the versification of some fundamental assumptions conveniently  made in the literature.
 The major contributions of this article thus include \vspace*{-0mm}
 \bi
 \item[(a)] introducing a novel partition of the parameter space and defining an original population version of the LTS for the first time;\vspace*{-2mm}
 \item[(b)] investigating  primary properties of the sample and population versions of the objective function for the LTS, obtaining original results;\vspace*{-2mm}
 \item[(c)] first time obtaining the influence function ($p\geq 2$) and Fisher consistency for the LTS; \vspace*{-2mm}
 \item[(d)] first time establishing the strong consistency of the sample LTS via a generalized Glivenko-Cantelli Theorem without artificial assumptions; and \vspace*{-2mm}
 \item[(e)]  first time employing a novel and concise approach based on the empirical process theory to establish asymptotic normality of the sample LTS. \vspace*{-2mm}
 \ei
\vs

The rest of the article is organized as follows. Section \ref{Sec.2} introduces first time the population version of LTS and addresses the properties of the LTS estimator in both empirical and population settings, including the global continuity and local differentiability and convexity of its objective function;
its influence function (in $p>2$ cases) and  Fisher consistency are established for the {first time}. Section \ref{Sec.3} establishes the strong consistency via a generalized Glivenko-Cantelli Theorem and  the
asymptotic normality of the estimator is re-established in a very different and concise approach (via stochastic  equicontinuity) rather than the previous ones in the literature.
Section \ref{Sec.4} addresses the asymptotic inference procedures based on the asymptotic normality and bootstrapping.
Concluding remarks in Section \ref{Sec.5} end the article. Major proofs are deferred to an Appendix.

\section {Definition and Properties of the LTS} \label{Sec.2}

\noin
\subsection{Definition} \label{Sec.2.1}
Denote by $F_{(\mb{x}^{\top}, y)}$ the joint distribution of $\mb{x}^{\top}$ and $y$ in model (\ref{model.eqn}). Throughout $F_{\mb{Z}}$ stands for the distribution function of the random vector $\mb{Z}^{\top}$.
For a given $\bs{\beta}\in \R^p$ and an $\alpha \in [1/2,c]$,  $1/2<c < 1$, let $q(\bs{\beta}, \alpha)=F^{-1}_W (\alpha)$ be the $\alpha$th quantile of  
 $F_W$ with $W:=W(\bs{\beta})=(y-\bs{w}^{\top}\bs{\beta})^2$,
where $\bs{w}^{\top}=(1,\bs{x}^{\top})$. The constant $c=1$ case is excluded to avoid  unbounded $q(\bs{\beta}, \alpha)$ and the LS  cases.
 Define an objective function
\be
O(F_{(\mb{x}^{\top}, y)}, \bs{\beta}, \alpha)=\int (y-\bs{w}^{\top}\bs{\beta})^2\mathds{1} \left((y-\bs{w}^{\top}\bs{\beta})^2\leq q(\bs{\beta}, \alpha)\right)dF_{(\mb{x}^{\top}, y)}, \label{objective.eqn}
\ee
and a regression functional
\be
\bs{\beta}_{lts}(F_{(\mb{x}^{\top}, y)},\alpha)=\arg\min_{\bs{\beta}\in \R^p}O(F_{(\mb{x}^{\top}, y)}, \bs{\beta}, \alpha), \label{lts-population.eqn}
\ee
 where $\mathds{1}(A)$ is the indicator of $A$ (i.e., it is one if A holds and zero otherwise).
Let $F^n_{(\mb{x}^{\top}, y)}$ be the sample version of the $F_{(\mb{x}^{\top}, y)}$ based on a sample $\mb{z}^{(n)}:=\{(\mb{x}^{\top}_i, y_i)^{\top}, i\in\{1,2,\cdots, n\}\}$. The $F^n_{(\mb{x}^{\top}, y)}$ and $\mb{z}^{(n)}$ will be used interchangeably. Using the $F^n_{(\mb{x}^{\top}, y)}$, one obtains the sample versions
\be
O(F^n_{(\mb{x}^{\top}, y)}, \bs{\beta}, \alpha)= \frac{1}{n}\sum_{i=1}^{\lfloor \alpha n \rfloor+1} r^2_{i:n}, \label{objective-0.eqn}
\ee
where $\lfloor x \rfloor$ is the floor function. Further
\be
\widehat{\bs{\beta}}^n_{lts}:=\bs{\beta}_{lts}(F^n_{(\mb{x}^{\top}, y)},\alpha) =\arg\min_{\bs{\beta}\in \R^p}O(F^n_{(\mb{x}^{\top}, y)}, \bs{\beta}, \alpha). \label{lts-1.eqn}
\ee
It is readily seen that the $\widehat{\bs{\beta}}^n_{lts}$ above is identical to the $\widehat{\bs{\beta}}_{lts}$ in (\ref{lts.eqn}) with $h=\lfloor \alpha n\rfloor+1$.
Henceforth we prefer to treat the $\widehat{\bs{\beta}}^n_{lts}$ rather than the $\widehat{\bs{\beta}}_{lts}$ in (\ref{lts.eqn}).\vs The first natural question is the existence of the minimizer in the right-hand side (RHS) of (\ref{lts-1.eqn}), or the existence of the $\widehat{\bs{\beta}}^n_{lts}$. Does it always exist? If it exists, will it be unique? Unique existence is a key precondition for the study of asymptotics of an estimator.
\vs
One might take the existence for granted since the objective is non-negative and has a finite infimum
which can be approximated by objective values of a sequence of $\bs{\beta}$s. There is a sub-sequence of $\bs{\beta}$s with their objective values converging to the infimum which is minimum
due the continuity of the objective function.  The sub-sequence of $\bs{\beta}$s converges to a point $\bs{\beta}^*$ which is the minimizer of the RHS.  There are multiple issues with the arguments above.
The existence and the convergence of the $\bs{\beta}$ sub-sequence (to a minimum) and continuity of objective function need to be proved. In the sequel, we take a different approach.

\noindent
\subsection{Properties in the empirical case} \label{Sec.2.2}
Write $O^n(\bs{\beta})$  and $\mathds{1}_i$
 for the $O(F^n_{(\mb{x}', y)},  \bs{\beta}, \alpha)$  and the $\mathds{1} \left(r^2_i\leq r^2_{h : n}\right)$, respectively.
  It is  seen that
\be
O^n(\bs{\beta})=\frac{1}{n}\sum_{i=1}^n r^2_i\mathds{1} \left(r^2_i\leq r^2_{h : n}\right)=\frac{1}{n}\sum_{i=1}^n r^2_i\mathds{1}_i,
\label{objective-1.eqn}
\ee
where $h=\lfloor \alpha n\rfloor+1$. The fraction $1/n$ will often be ignored in the following discussion.

\vs
\noindent
\tb{Existence and uniqueness}\vs
\noin
\tb{Partition parameter space}~~ 
 For a given sample $\bs{z}^{(n)}$, 
 an $\alpha$ (or $h$), and any $\bs{\beta}^1\in\R^p$, 
 let  $r^2_{i:n}(\bs{\beta}^1)=r^2_{k_i}(\bs{\beta}^1)$ for an integer $k_i$.  Note that $r_j$ and $k_i$ depend on $\bs{\beta}^1$, i.e., $r_j:=r_j(\bs{\beta}^1)$, $k_i:=k_i(\bs{\beta}^1)$. 
Obviously 
 $r^2_{k_1}\leq r^2_{k_2}\leq \cdots\leq r^2_{k_n}$. 
 Call $\{k_i, 1\leq i\leq h\}$ $\bs{\beta}^1$-$h$-integer set. 
   If 
    {$r_i^2\not =r^2_j$ for any distinct $i$ and $j$},     
   then the $h$-integer set is unique. Hereafter we assume  
       \tb{(A0):}\\ \tb{$W$ has a density for any given $\bs{\beta}$}.
  Then almost surely (a.s.), the $h$-integer set is unique. 
  \vs
   Consider the unique cases. There are other $\bs{\beta}$s in $\R^p$ that share the same $h$-integer set as that of the $\bs{\beta}^1$. Denote the set of such points that have the same h integers as $\bs{\beta}^1$ by
 \be
 S_{\bs{\beta}^1}:=\Big\{\bs{\beta}\in \R^p: k_i(\bs{\beta})=k_i(\bs{\beta^1})=k_i,~ i\in\{1, 2, \cdots, h\}~~
 \{k_i, 1\leq i\leq h\}  ~\mbox{is unique.}\Big\} 
 \label{element-set.eqn}
 \ee
\vs
 If \tb{(A0)} holds, then $S_{\bs{\beta}^1} \neq \emptyset$ (a.s.). If it is
 $\R^p$, then we have a trivial case (see Remarks 2.1 below).  Otherwise, there are only finitely many such sets (for a fixed $n$) that partition  
 $\R^p$. Let $\cup_{l=1}^L \overline{S}_{\bs{\beta}^l}=\R^p$, where 
  ${S}_{\bs{\beta}^l}$s are defined similarly to (\ref{element-set.eqn}) and 
   are disjoint for different $l$,  $1\leq l\leq L:= {n \choose h}$, and $\overline {A} $ is the closure of the set $A$. Write $\mb{X}_n=(\bs{w}_1, \cdots, \bs{w}_n)'$, an $n\times p$ matrix. Assume \tb{(A1):} \tb{$\mb{X}_n$ and any its h rows have a full rank $p$}. As in the 
 \textsf{R}
  package ltsReg (by Valentin Todorov), hereafter we assume that $p<n/2$.
 \vs
\noindent
\tb{Lemma 2.1}
Assume that \tb{(A0)} and \tb{(A1)} hold, then
\bi
\item[(i)] \tb{(a)} For any $l$ ($1\leq l\leq L$),  $r^2_{k_1(\bs{\beta}^l)}< r^2_{k_2(\bs{\beta}^l)}<\cdots <r^2_{k_h(\bs{\beta}^l)}$ over  $S_{\bs{\beta}^l}$. \\
\tb{(b)} For 
any $\bs{\eta}\in S_{\bs{\beta}^l}$, there exists
     an open ball $B(\bs{\eta}, \delta)$ centered at $\bs{\eta}$ with a radius $\delta>0$ such that for any $\bs{\beta}\in B(\bs{\eta}, \delta)$ 
    \be O^n(\bs{\beta})=\frac{1}{n}\sum_{i=1}^hr^2_{k_i (\bs{\beta})} (\bs{\beta})=\frac{1}{n}\sum_{i=1}^hr^2_{k_i (\bs{\eta})} (\bs{\beta})=\frac{1}{n}\sum_{i=1}^hr^2_{k_i (\bs{\beta}^l)} (\bs{\beta}), ~(a.s.)\label{expression-4-O-n.eqn}
    \ee
\item[(ii)] The graph of $O^n(\bs{\beta})$ over $\bs{\beta}\in \R^p$ is composed of the $L$ closures of
 graphs of the quadratic function of $\bs{\beta}$: $
\frac{1}{n}\sum_{i=1}^hr^2_{k_i (\bs{\beta}^l)} (\bs{\beta})$ for $ O^n(\bs{\beta}^l) $ and
 any $l$ ($1\leq l\leq L$), joined together.
\item[(iii)] $O^n(\bs{\beta})$ is continuous in $\bs{\beta}\in \R^p$.
\item[(iv)] $O^n(\bs{\beta})$ is differentiable and strictly convex over each $S_{\bs{\beta}^l}$ for any $1\leq l\leq L$. 
\ei

\noindent
\tb{Proof:} see the Appendix. \hfill $\square$  
\vs\noindent
\tb{Remarks 2.1}\vs
\tb{(a)} If $S_{\bs{\beta}^0}=\R^p$, then $O^n(\bs{\beta})$ is a twice differentiable and strictly convex quadratic function of $\bs{\beta}$, the existence and the uniqueness of $\widehat{\bs{\beta}}^n_{lts}$ are trivial as long as $\mb{X}_n$ has a full rank.
\vs
 \tb{(b)}  Replacing $(1,\bs{x}_i)'\bs{\beta}$  by a nonlinear $h(\bs{x}_i, \bs{\beta})$, 
  C04 and C05 
also addressed the continuity and differentiability of the objective function of the LTS.
However, they assumed that (i) $F_W$ is twice differentiable around the points corresponding to the square roots
of the $\alpha$-quantiles of $W$, (ii) $h(\bs{x}_i, \bs{\beta})$ is continuous over parameter space $B$, (iii) $h(\bs{x}_i, \bs{\beta})$ is  twice differentiable in $\bs{\beta}$ for $\bs{\beta}\in B(\bs{\beta}_0, \delta)$ a.s., and (iv) $\partial h(\bs{x}_i,\bs{\beta})/\partial\bs{\beta}$ is continuous in $\bs{\beta}$,
All  assumptions were never verified in C04 and C05 though, 
 however they are not required in Lemma 2.1. \vs


 \tb{(c)} Continuity and differentiability inferred just based on $O^n(\bs{\beta})$ being the sum of $h$ continuous and differentiable functions (squares of residual) without (i) or (\ref{expression-4-O-n.eqn}) might be inadequate (flawed). In general, $O^n(\bs{\beta})$ is not differentiable nor convex in $\bs{\beta}$ {globally}.
{\hfill \ensuremath{\Box}} 
\vs
\vs
Let $\mb{y}_n:=(y_1,\cdots, y_n)^{\top}$,  $\bs{M}_n:=\bs{M}(\mb{y}_n, \mb{X}_n, \bs{\beta}, \alpha)=\sum_{i=1}^n\bs{w}_i\bs{w}^{\top}_i\mathds{1}_i=\sum_{i=1}^h \bs{w}_{k_i(\bs{\beta})}\bs{w}^{\top}_{k_i(\bs{\beta})}$. Note that $\mathds{1}_i$ depends on $\bs{\beta}$. 
\vs
\noindent
\tb{Theorem 2.1}  Assume that \tb{(A0)} and \tb{(A1)} hold. Then
\bi
\item[(i)] $\widehat{\bs{\beta}}^n_{lts}$ exists and is the local minimum of $O^n(\bs{\beta})$ over $S_{\bs{\beta}^{l_0}}$ for some $l_0$ ($1\leq l_0\leq L$).

\item[(ii)] 
Over $S_{\bs{\beta}^{l_0}}$,  $\widehat{\bs{\beta}}^n_{lts}$  is the solution of the system of equations
\be
 \sum_{i=1}^{n}(y_i-\bs{w}^{\top}_i\bs{\beta})\bs{w}_i\mathds{1}_i 
  =\mb{0}, \label{estimation.eqn}
\ee
\item[(iii)]
Over $S_{\bs{\beta}^{l_0}}$, the unique solution is 
\be
\widehat{\bs{\beta}}^n_{lts}=\bs{M}_n(\mb{y}_n, \mb{X}_n, \widehat{\bs{\beta}}^n_{lts}, \alpha)^{-1}\sum_{i=1}^hy_{k_i(\bs{\beta}^{l_0})}\bs{w}_{k_i(\bs{\beta}^{l_0})}.
\label{lts-iterative.eqn}
\ee
\vspace*{-6mm}
\ei

\noindent
\tb{Proof:} The given conditions and Lemma 2.1 allow one to focus on a piece $S_{\bs{\beta}^l}$, $1\leq l\leq L$, all results follows in a straightforward fashion, for more details, see the Appendix. \hfill \pend
\vs
\noindent
\tb{Remarks 2.2}
\vs
\tb{(a)} Unique existence,  which is often implicitly assumed or ignored in the literature,  is central for the discussion of asymptotics of $\widehat{\bs{\beta}}^n_{lts}$.
 Existence of $ \widehat{\bs{\beta}}^n_{lts}$ could also be established under the assumption that there are no $\lfloor (n+1)/2\rfloor$ sample points of $\bs{z}^{(n)}$ contained in any $(p-1)$-dimensional hyperplane, similarly to that of Theorem 2.2 for LST in Zuo and Zuo \cite{ZZ23}. It is established here without such assumption nevertheless.
\vs
\tb{(b)} A sufficient condition for the invertibility of $\bs{M}_n$ is that any $h$ rows of $\bs{X}_n$ form a full rank sub-matrix. The latter is true if  \tb{(A1)} holds.
\vs
\tb{(c)} V\'{i}\v{s}ek (2006a)  
  also addressed  the existence of $ \widehat{\bs{\beta}}^n_{lts}$ (Assertion 1) for non-random covariates (carriers) satisfying many demanding assumptions ($\mathcal{A, B}$).
   The uniqueness was left unaddressed though.
\hfill \pend

\vs
\noin
\subsection{Properties in the population case}
\vs
The best breakdown point of the LTS (see p. 132 of RL87) 
reflects its global robustness.
We now examine its local robustness via the influence function to depict its complete robust picture.
\vs
\noin
\tb{Definition of influence function} 
\vs
For a distribution $F$ on $\R^{p}$
and an $\varepsilon \in (0, 1/2)$, the version of $F$ contaminated by an $\varepsilon$ amount of an arbitrary
distribution $G$ on $\R^{p}$  is denoted by $F(\varepsilon, G) = (1-\varepsilon)F + \varepsilon G$ (an $\varepsilon$ amount deviation from
the assumed $F$). $F(\varepsilon,G)$ is actually a convex contamination of $F$. There
are other types of contamination such as contamination by total variation or Hellinger distances. 
 We cite the definition of influence function  from Hampel et al. \cite{HRRS86}.
\vs
\noindent
\tb{Definition 2.2}~\cite{HRRS86} 
The \emph{influence function} (IF) of a functional $\mb{t}$ at a given point $\bs{x} \in \R^{p}$ for a given $F$ is defined
as
\be
\text{IF}(\bs{x}; \mb{t},F) = \lim_ {\ep\to 0^+} \frac{\mb{t} (F (\ep, \delta_{\bs{x}})) - \mb{t} (F )}{\ep},\label{if.eqn} \ee
where $\delta_{\bs{x}}$ is the point-mass probability measure at $\bs{x} \in \R^{p}$.
\vs
The function $\mbox{IF}(\bs{x}; \mb{t}, F)$ describes the relative effect (influence) on $\mb{t}$ of an infinitesimal
point-mass contamination at $\bs{x}$ and measures the local robustness of $\mb{t}$.\vs

To establish the IF for the functional $\bs{\beta}_{lts}(F_{(\bs{x}^{\top}, y)}, \alpha)$, we need to first show its existence and uniqueness with or without point-mass contamination. To that end, write
$$F_{\varepsilon}(\mb{z}):=F(\varepsilon, \delta_{\mb{z}})= (1-\varepsilon) F_{(\bs{x^{\top}}, y)}+\varepsilon \delta_{\mb{z}},$$
with $\mb{u}=(\mb{s}^{\top}, t)^{\top}\in \R^{p}$, $\mb{s} \in \R^{p-1}$ , $t\in \R^1$ being the corresponding random vector (i.e. $F_{\mb{u}}=F_{\varepsilon}(\mb{z})=F(\varepsilon, \delta_{\mb{z}})$).
The versions of (\ref{objective.eqn}) and (\ref{lts-population.eqn}) at the contaminated $F(\varepsilon, \delta_{\mb{z}})$ are respectively
\be
O(F_{\varepsilon}(\mb{z}), \bs{\beta}, \alpha)=\int(t-\bs{v}^{\top}\bs{\beta})^2\mathds{1} \left((t-\bs{v}^{\top}\bs{\beta})^2\leq q_{\varepsilon}(\mb{z}, \bs{\beta},\alpha)\right)dF_{\mb{u}}(\mb{s}, t), \label{objective-contaminated.eqn}
\ee
 with $q_{\varepsilon}(\mb{z},\bs{\beta},\alpha)$ being the $\alpha$th quantile of the distribution function of $(t-\bs{v}^{\top}\bs{\beta})^2$, $\bs{v}=(1,\bs{s}^{\top})^{\top}$, and
 \be
 \bs{\beta}_{lts}(F_{\varepsilon}(\mb{z}), \alpha)=\arg\min_{\bs{\beta}\in \R^p}O(F_{\varepsilon}(\mb{z}), \bs{\beta}, \alpha). \label{lst-contaminated.eqn}
 \ee
 \vs
 For $\bs{\beta}_{lts}(F_{(\bs{x}^{\top}, y)}, \alpha)$ defined in (\ref{lts-population.eqn}) and $\bs{\beta}_{lts}(F_{\varepsilon}(\mb{z}), \alpha)$ above, we have an analogous result to Theorem 2.1. (Assume that the counterpart of model (\ref{model.eqn}) is $y=(1, \bs{x}^{\top})\bs{\beta}_0+e=\bs{w}^{\top}\bs{\beta}_0+e$). Before we derive the influence function, we need to establish existence and uniqueness.
 \vs
 \vs
 \noindent
 \tb{Existence and Uniqueness}\vs
\noindent
Write $O(\bs{\beta})$ for $O(F_{(\bs{x}^{\top}, y)}, \bs{\beta}, \alpha)$ in (\ref{objective.eqn}). To have a counterpart of Lemma 2.1, we need\\[2ex]
\tb{(A2)} $W$ has a positive density around a small neighborhood of $q(\bs{\beta}, \alpha))$ for the given $\alpha$, $\bs{\beta}$.
  
\vs
\noindent
\tb{Lemma 2.2} Assume \tb{(A2)} holds 
and $E(\bs{w}\bs{w}^{\top})$ exists. Then %
\vs
\tb{(i)} $O(\bs{\beta})$ is continuous in $\bs{\beta} \in \R^p$;
\vs
\tb{(ii)} $O(\bs{\beta})$ is twice differentiable  in $\bs{\beta} \in \R^p$, \[{\partial^2 O(\bs{\beta})}\big/{\partial \bs{\beta}^2}=2E(\bs{w}\bs{w}^{\top} \mathds{1}((y-\bs{w}^{\top}\bs{\beta})^2\leq q(\bs{\beta}, \alpha));\]
\vs
\tb{(iii)}  $O(\bs{\beta})$ is strictly convex in $\bs{\beta} \in \R^p$. 
 \vs \noindent
\tb{Proof:} The boundedness of the integrand in (\ref{objective.eqn}), given conditions and Lebesgue dominated convergence theorem leads to the desired results, for details, see the Appendix. \hfill \pend
 \vs
Note that \tb{(ii)} and \tb{(iii)} above are global in $\bs{\beta}$,  stronger than the empirical counterparts, 
 all are attributed to boundary of $S_{\bs{\beta}^l}$ issue. We now treat
the existence and uniqueness of $\bs{\beta}_{lts}$, which is central for the asymptotics study. \vs
 \noin
 \tb{Theorem 2.3} Assume  \tb{(A2)} holds  and $E(\bs{w}\bs{w}^{\top})$ exists,
 $q_{\varepsilon}(\mb{z}, \bs{\beta},\alpha)$ is continuous in $\bs{\beta}$, and
 $P((t-\bs{v}^{\top}\bs{\beta})^2=q_{\varepsilon}(\mb{z}, \bs{\beta},\alpha))=0$ for any $\bs{\beta}\in \R^p$ and the given $\alpha$. Then 
 \vs
 \tb{(i)} $\bs{\beta}_{lts}(F_{(\bs{x}^{\top}, y)}, \alpha)$  and $\bs{\beta}_{lts}(F_{\varepsilon}(\mb{z}), \alpha)$  exist.
 \vs
 \tb{(ii)} Furthermore, they are the solution of system of equations,  respectively
 \begin{align}
 \int(y-\bs{w}^{\top}\bs{\beta})\bs{w}\mathds{1}\left((y-\bs{w}^{\top}\bs{\beta})^2\leq q(\bs{\beta}, \alpha) \right) dF_{(\bs{x}^{\top}, y)}(\bs{x}, y)&=\mb{0},  \label{lst-estimation.eqn} \\
 \int(t-\bs{v}^{\top}\bs{\beta})\bs{v}\mathds{1}\left((t-\bs{v}^{\top}\bs{\beta})^2\leq q_{\varepsilon}(\mb{z},\bs{\beta},\alpha) \right) dF_{\mb{u}}(\mb{s}, t)&=\mb{0}. \label{lts-contamination-estimation.eqn}
 \end{align}

 \tb{(iii)} $\bs{\beta}_{lts}(F_{(\bs{x}^{\top}, y)}, \alpha)$  and $\bs{\beta}_{lts}(F_{\varepsilon}(\mb{z}), \alpha)$ are unique
  provided that
 \begin{align}
  \int \bs{w}\bs{w}^{\top}\mathds{1}\left((y-\bs{w}^{\top}\bs{\beta})^2\leq q(\bs{\beta}, \alpha) \right)dF_{(\bs{x}^{\top}, y)}(\bs{x}, y)&, \label{uniqueness-lts.eqn}\\
  \int \bs{v}\bs{v}^{\top}\mathds{1}\left((t-\bs{v}^{\top}\bs{\beta})^2\leq q_{\varepsilon}(\mb{z},\bs{\beta},\alpha) \right) dF_{\mb{u}}(\mb{s}, t)&
  \end{align}
  {are invertible for $\bs{\beta}$ in a small neighborhood of $\bs{\beta}_{lts}(F_{(\bs{x}^{\top}, y)}, \alpha)$  and $\bs{\beta}_{lts}(F_{\varepsilon}(\mb{z}), \alpha)$ respectively.}
 \vs
 \noindent
\tb{Proof:} In light of Lemma 2.2, the proof is straightforward, see the Appendix. \hfill \pend
 \vs
 Continuity of $q_{\varepsilon}(\mb{z}, \bs{\beta},\alpha)$ in $\bs{\beta}$ is necessary for the
  differentiability of $O(F_{\varepsilon}(\mb{z}), \bs{\beta}, \alpha)$. In the non-contaminated case, continuity of $q(\bs{\beta}, \alpha)$ is guaranteed by \tb{(A2)}.\vs
 Does the population version of the LTS, $\bs{\beta}_{lts}$, defined in (\ref{lts-population.eqn}), have something to do with $\bs{\beta}_0$? It turns out
 under some conditions, they are identical, which is called Fisher consistency.
 \vs
 \noindent
 \tb{Fisher consistency} \vs \noin
 \tb{Theorem 2.4} Assume 
  \tb{(A2)} holds and $E(\bs{w}\bs{w}^{\top})$ exists,
  then $\bs{\beta}_{lts}(F_{(\bs{x}^{\top}, y)}, \alpha)=\bs{\beta}_0$  provided that
(i) $E_{(\bs{x}^{\top}, y)}\left(\bs{w}\bs{w}^{\top} \mathds{1}(r(\bs{\beta})^2\leq F^{-1}_{r(\bs{\beta})^2}(\alpha))\right)$ is invertible,
and\\
 (ii) $E_{(\bs{x}^{\top}, y)}(e\bs{w}\mathds{1}(e^2\leq F^{-1}_{e^2}(\alpha))=\mb{0}$,
 where $r(\bs{\beta})=y-\bs{w}^{\top}\bs{\beta}$. 
\vs \noindent
\tb{Proof:} Theorem 2.3 leads directly to the desired result, see the Appendix. \hfill \pend
 \vs
 \vs
 \noindent
\tb{Influence function}\vs
\noin
 \tb{Theorem 2.5} Assume that the assumptions in Theorem 2.3 hold. Set $\bs{\beta}_{lts}:=\bs{\beta}_{lts}(F_{(\bs{x}^{\top}, y)}, \alpha)$. Then for any $\mb{z}_0:=(\mb{s}^{\top}_0, t_0)^{\top} \in \R^{p}$, we have that
 \[
 \mbox{IF}(\mb{z}_0; \bs{\beta}_{lts}, F_{(\bs{x}^{\top}, y)})=\left\{
 \begin{array}{ll}
 \mb{0}, & \mbox{if $(t_0-\bs{v}^{\top}_0\bs{\beta}_{lts})^2>q(\bs{\beta}_{lts}, \alpha)$,}\\[1ex]
 M^{-1}(t_0-\bs{v}^{\top}_0 \bs{\beta}_{lts})\bs{v}_0,&    \mbox{otherwise},
 \end{array}
 \right.
 \]
 provided that $\bs{M}=E_{(\bs{x}^{\top}, y)}\Big(\bs{w}\bs{w}^{\top}\mathds{1}\left(r(\bs{\beta}_{lts})^2\leq q(\bs{\beta}_{lts},\alpha) \right)\Big)$ is invertible, where $\bs{v}_0=(1,\bs{s}^{\top}_0)^{\top}$.
 \vs

 \noindent
\tb{Proof:} The connection to the derivative of a functional is the key, see the Appendix. \hfill \pend
 \vs
 \noindent
 \tb{Remarks 2.3}
 \vs
 \tb{(a)} When $p=1$, the problem in  our model (\ref{model.eqn}) becomes a location problem (see p. 158 of RL87) 
 and the IF of the LTS estimation functional has been given on p. 191 of  RL87. 
 In the location setting, Tableman \cite{T94} also studied the IF of the LTS. When $p=2$, namely in the simple regression case, \"{O}llerer et al. (2015) 
 studied IF of the sparse-LTS functional under the assumption
 that $\bs{x}$ and ${e}$ are independent and normally distributed. Under stringent assumptions on the error terms $e_i$ and on  $\bs{x}$, Ma\v{s}\'{i}\v{c}ek (2004) 
 also addressed the IF of LTS for any $p$, but the point mass at $(\bs{x}^{\top}, z)$ with $z$ being the error term, an unusual contaminating point.
 The above result is much more general and valid for any $p\geq 1$, $\bs{x}'$, and $e$.
 \vs
 \tb{(b)} The influence function of $\bs{\beta}_{lts}$ remains
bounded if the contaminating point $(\bs{s}'_0, t_0)$ does not follow the model (i.e., its residual is extremely large), in particular for bad leverage points and vertical outliers. This shows the
good robust properties of the LTS.
 \vs
\tb{(c)} The influence function of $\bs{\beta}_{lts}$, unfortunately, might be unbounded (in $p > 1$ case),
sharing the drawback of the sparse-LTS (in the p = 2 case). The latter was shown in \"{O}llerer et al. (2015). 
Trimming based on the residuals (or squared residuals) will have this type of
drawback since the term $\bs{w}^{\top}\bs{\beta}$ can be bounded, but $\|\bs{x}\|$ might not.
 \hfill \pend

\section{Asymptotic properties} \label{Sec.3}

V\'{i}\v{s}{e}k (2006a, b, c) (V06a, b, c) 
 rigorously addressed the consistency, root-n consistency, and normality of the LTS under a restrictive setting ($\bs{x}_i$s are non-random covariates) plus 
  many assumptions on $\bs{x}_i$s and on the distribution of $e_i$ in three lengthy series papers. \vs\v{C}\'i\v{z}ek (2004, 2005) (C04, C05) 
   also addressed asymptotic properties of an extended LTS under $\bs{\beta}$-mixing conditions for $\bs{x}_i$ with nonlinear regression function $h(\bs{x}_i,\bs{\beta})$, a 
    more general setting, but with the price of numerous artificial assumptions (H1, H2, H3, H4, H5, H6; D1,D2; I1,I2) that are never verified in any concrete example 
nor 
for the special case of linear LTS. That is, C04 and C05 cannot trivially be applied for  the linear LTS in (\ref{lts.eqn}).\vs

 Strong consistency has been addressed in CSZ97 in nonlinear setting with  key assumptions 
 that are never verified even for the special linear LTS.
We now rigorously establish strong consistency without any artificial assumptions.
\vs
\noindent
\subsection{Strong consistency}
\vs
Following the notations of Pollard (1984) (P84), 
write
\begin{align*}
O(\bs{\beta}, P):&=O(F_{(\bs{x}^{\top}, y)}, \bs{\beta}, \alpha)=P [(y-\bs{w}^{\top}\bs{\beta})^2\mathds{1}(r(\bs{\beta})^2\leq F^{-1}_{r(\bs{\beta})^2}(\alpha))]=Pf,\\[1ex]
 O(\bs{\beta}, P_n):&=O(F^n_{(\bs{x}^{\top}, y)}, \bs{\beta}, \alpha)
 =\frac{1}{n}\sum_{i=1}^n r^2_i\mathds{1} \left(r^2_i\leq (r)^2_{h : n}\right)= P_nf,
\end{align*}
where $f:=f(\bs{x}, y, \bs{\beta}, \alpha)=(y-\bs{w}^{\top}\bs{\beta})^2\mathds{1}(r(\bs{\beta})^2\leq F^{-1}_{r(\bs{\beta})^2}(\alpha))$, $\alpha$ and $h=\lfloor \alpha n \rfloor+1$ are fixed. 
\vs
Under corresponding assumptions in Theorems 2.1 and 2.3, $\widehat{\bs{\beta}}^n_{lts}$ and $\bs{\beta}_{lts}$ are unique minimizers of $O(\bs{\beta}, P_n)$ and $O(\bs{\beta}, P)$,
respectively.\vs
To show that $\widehat{\bs{\beta}}^n_{lts}$ converges to $\bs{\beta}_{lts}$ a.s., one can take the approach given in Section 4.2 of Zuo \cite{ZZ23}. But here we take a different and more direct approach.\vs
To show that $\widehat{\bs{\beta}}^n_{lts}$ converges to $\bs{\beta}_{lts}$ a.s., it will suffice to prove that $O(\widehat{\bs{\beta}}^n_{lts}, P) \to O(\bs{\beta}_{lts}, P)$ a.s., because $O(\bs{\beta}, P)$ is bounded away from $O(\bs{\beta}_{lts}, P)$ outside each neighborhood of $\bs{\beta}_{lts}$ in light of continuity and compactness (also see Lemma 4.3 of Zuo \cite{ZZ23}).
Let $\Theta$ be a closed ball centered at $\bs{\beta}_{lts}$ with a radius $r>0$. Define a class of functions
\[
\mathscr{F} (\bs{\beta}, \alpha)=\left\{f(\bs{x}, y, \bs{\beta}, \alpha)= (y-\bs{w}^{\top}\bs{\beta})^2\mathds{1}(r(\bs{\beta})^2\leq F^{-1}_{r(\bs{\beta})^2}(\alpha)): \bs{\beta} \in \Theta, \alpha \in [1/2, c]\right\}
\]
\vs
\noin
If we prove uniform almost sure convergence of $P_n$ to $P$ over $\mathscr{F}$ (see Lemma 3.1 below), then we can deduce that $O(\widehat{\bs{\beta}}^n_{lts}, P) \to O(\bs{\beta}_{lts}, P)$  a.s.  from
\begin{align*}
O(\widehat{\bs{\beta}}^n_{lts}, P_n)-O(\widehat{\bs{\beta}}^n_{lts}, P) &\to 0 ~~(\mbox{in light of Lemma 3.1), ~~and}\\[1ex]
O(\widehat{\bs{\beta}}^n_{lts}, P_n)\leq O(\bs{\beta}_{lts}, P_n)& \to O(\bs{\beta}_{lts}, P)\leq O(\widehat{\bs{\beta}}^n_{lts}, P). 
\end{align*}

The above discussions and arguments have led to
\vs
\noindent
\tb{Theorem 3.1}. Under corresponding assumptions in Theorems 2.1 and 2.3 for uniqueness of $\widehat{\bs{\beta}}^n_{lts}$ and $\bs{\beta}_{lts}$ respectively, we have
$\widehat{\bs{\beta}}^n_{lts}$ converges a.s. to $\bs{\beta}_{lts}$ (i.e. $\widehat{\bs{\beta}}^n_{lts}-\bs{\beta}_{lts}=o(1)$, a.s.).
\vs
The above is based on the following generalized Glivenko-Cantelli Theorem.
\vs
\noindent
\tb{Lemma 3.1}. $\sup_{f\in \mathscr{F}}|P_n f-P f| \to 0$ a.s., provided that \tb{(A2)} holds.
\vs
\noindent
\vs \noindent
\tb{Proof:} Verifying two requirements in the Theorem 24 in II.5 of P84 leads to the result. Showing the covering number for functions in $\mathscr{F}$ is bounded is challenging. Essentially
one needs to show that the graphs of functions in $\mathscr{F}$ form a VC class of sets (this was avoided in the literature, e.g., in C04, C05, CSZ97, and V06a, b, c). For details,
 see the Appendix. \hfill \pend
\vs
\noindent
\subsection{Root-n consistency and asymptotic normality} \vs
Instead of treating the root-n consistency separately as in V06b c, we will establish asymptotic normality of $\widehat{\bs{\beta}}^n_{lts}$ directly via stochastic equicontinuity (see p. 139 of P84 or the supplementary of Zuo (2020)). 
\vs
\emph{Stochastic equicontinuity} refers to a sequence of
stochastic processes $\{Z_n(t): t \in T\}$ whose shared index set $T$ comes equipped
with a semi-metric $d(\cdot, \cdot)$. 
\vs
\noin
\tb{Definition 3.1} [IIV. 1, Def. 2 of {P84}].  Call ${Z_n}$ stochastically equicontinuous at $t_0$  if for each $\eta > 0$
and $\epsilon > 0$ there exists a neighborhood $U$ of $t_0$ for which
\be
\limsup P\left(\sup_{U} |Z_n(t) - Z_n(t_0) | > \eta\right) < \epsilon.  \label{se.eqn}
\ee
\vs
It is readily seen (see p.140 of P84) that
if ${\tau_n}$ is a sequence of random elements of $T$ that converges in probability
to $t_0$, then
\be
Z_n(\tau_n)-Z_n(t_0)\to 0\mbox{~in probability,}
\ee
because, with probability tending to one, $\tau_n$ will belong to each $U$.
The form above will be easier to apply, especially when behavior of a particular ${\tau_n}$ sequence is under investigation.\vs
 Suppose $\mathscr{F} = \{ f(\cdot, t): t\in T\}$, with $T$ a subset of $\R^k$, is a collection of
real, P-integrable functions on the set $S$ where $P$ (probability measure) lives.
Denote by $P_n$ the
empirical measure formed from $n$ independent observations on $P$, and define the empirical process $E_n$ as the signed measure $n^{1/2}(P_n - P)$. Define
\begin{align*}
F(t) &= P f(\cdot, t),\\
F_n(t) &= P_n f(\cdot, t).
\end{align*}
Suppose $f(\cdot, t)$ has a linear approximation near the $t_0$ at which $F(\cdot)$ takes
on its minimum value:
\be
 f(\cdot, t) = f(\cdot, t_0) + (t - t_0)^{\top}\nabla  f(\cdot, t) + |t - t_0|r(\cdot, t). \label{taylor.eqn}
\ee
For completeness set $r(\cdot, t_0) = 0$, where $\nabla$ (differential operator) is a vector of $k$ real functions on
$S$. We cite theorem 5 in VII.1 of {P84} (p. 141) for the asymptotic normality of $\tau_n$.
\vs
\noindent
\tb{Lemma 3.2} [{P84}].  Suppose $\{\tau_n\}$ is a sequence of random vectors converging in
probability to the value $t_0$ at which $F(\cdot)$ has its minimum. Define $r(\cdot, t)$ and the
vector of functions $\nabla(\cdot, t)$ by (\ref{taylor.eqn}). If
\bi
\item[(i)] $t_0$ is an interior point of the parameter set $T$; \vspace*{-2mm}
\item[(ii)] $F(\cdot)$ has a non-singular second derivative matrix $V$ at $t_0$;\vspace*{-2mm}
\item[(iii)] $F_n(\tau_n) = o_p(n^{-1}) + \inf_{t}F_n(t)$;\vspace*{-2mm}
\item[(iv)] the components of $\nabla f(\cdot, t)$ all belong to $\mathscr{L}^2(P)$;\vspace*{-2mm}
\item[(v)] the sequence $\{E_{n}r(\cdot, t)\}$ is stochastically equicontinuous at $t_0$;\vspace*{-2mm}
\ei
then
\[n^{1/2}(\tau_n - t_0) \stackrel{d} \longrightarrow  {{N}}(\bs{0}, V^{-1}[P(\nabla\nabla^{\top}) - (P\nabla)(P\nabla)^{\top}]V^{-1}).
\]
\vs
By Theorems 2.1 and 3.1, assume, without loss of generality (w.l.o.g.), that  $\widehat{\bs{\beta}}^n_{lts}$ and $\bs{\beta}_{lts}$ belong to a ball centered at $\bs{\beta}_{lts}$ with a large enough radius $r_0$, $B(\bs{\beta}_{lts}, r_0)$, and that $\Theta=B(\bs{\beta}_{lts}, r_0)$ is our parameter space of $\bs{\beta}$ hereafter.
In order to apply the Lemma, we first realize that in our case, $\widehat{\bs{\beta}}^n_{lts}$ and $\bs{\beta}_{lts}$ correspond to $\tau_n$ and $t_0$ (assume, w.l.o.g. that $\bs{\beta}_{lts}=\mb{0}$ in light of regression equivariance, see Section 4); $\bs{\beta}$ and $\Theta$ correspond to $t$ and $T$; $f(\cdot, t):= f(\bs{x}^{\top}, y,  \bs{\beta}, \alpha)=(y-\bs{w}^{\top}\bs{\beta})^2\mathds{1}(r(\bs{\beta})^2\leq F^{-1}_{r(\bs{\beta})^2}(\alpha))$.
In our case (see the proof (ii) of Lemma 2.2 in the Appendix),
\[
\nabla f(\cdot, t):=
\nabla f(\bs{x}, y, \bs{\beta}, \alpha)=\frac{\partial}{\partial \bs{\beta}} f(\bs{x}, y, \bs{\beta}, \alpha)=-2(y-\bs{w}^{\top}\bs{\beta})\bs{w}\mathds{1}\left(r(\bs{\beta})^2\leq F^{-1}_{r(\bs{\beta})^2}(\alpha)\right).
\]
We will have to assume that $P(\nabla^2_i)=P(4(y-\bs{w}^{\top}\bs{\beta})^2{w}^2_i\mathds{1}(r(\bs{\beta})^2\leq F^{-1}_{r(\bs{\beta})^2}(\alpha)))$ exists to meet (iv) of the lemma, where $i\in \{1,\cdots, p\}$ and $\bs{w}^{\top}=(1,\bs{x}^{\top})=(1, x_1, \cdots, x_{p-1})$. It is readily seen that a sufficient condition for this assumption to hold is the existence of $P(x^2_i)$. In our case,
$V= 2P(\bs{w}\bs{w}^{\top}\mathds{1}(r(\bs{\beta})^2\leq F^{-1}_{r(\bs{\beta})^2}(\alpha)))$, and we will have to assume that it is invertible when $\bs{\beta}$ is replaced by $\bs{\beta}_{lts}$ (this is covered by (\ref{uniqueness-lts.eqn}))  to meet (ii) of the Lemma. In our case,
\[r(\cdot, t)=\left(\frac{\bs{\beta}^{\top}}{\|\bs{\beta}\|}V/2 \frac{\bs{\beta}}{\|\bs{\beta}\|} \right)\|\bs{\beta}\|.\]
We will assume that $\lambda_{min}$ and $\lambda_{max}$ are the minimum and maximum eigenvalues of positive semidefinite matrix $V$ over all $\bs{\beta}\in \Theta$ and $\alpha \in[1/2, c]$.
\vs
\noindent
\tb{Theorem 3.2} Assume that\vspace*{-2mm}
\bi
\item[(i)] the uniqueness assumptions for $\widehat{\bs{\beta}}^n_{lts}$ and $\bs{\beta}_{lts}$ in Theorems 2.1 and 2.3 hold respectively;\vspace*{-2mm}
\item[(ii)] $P({x^2_i})$ exists with $\bs{x}=(x_1, \cdots, x_{p-1})$; \vspace*{-2mm}
\ei
then
\[n^{1/2}(\widehat{\bs{\beta}}^n_{lts} - \bs{\beta}_{lts}) \stackrel{d} \longrightarrow  {{N}}(\bs{0}, V^{-1}[P(\nabla\nabla^{\top}) - (P\nabla)(P\nabla)^{\top}]V^{-1}),
\]
where $\bs{\beta}$ in $V$ and $\nabla$ is replaced by $\bs{\beta}_{lts}$ (which could be assumed to be zero, or is $\bs{\beta}_0$ under Theorem 2.4).
\vs \noindent
\tb{Proof:} The key to apply Lemma  is to verify (v), for details  see the Appendix. \hfill \pend
\vs
\noindent
\tb{Remarks 3.1}
\vs
\tb{(a)} In the case of $p=1$, that is in the location case, the asymptotic normality of the LTS has been studied in Butler (1982), RL87, 
  Bednarski and Clarke (1993), 
  and Ma\v{s}\'{i}\v{c}ek (2004). 
  \vs

\tb{(b)} H\"{o}ssjer (1994), 
under the rank-based optimization framework and stringent assumptions on error term $e_i$ (even density that is strictly decreasing for positive value, bounded absolute first moment)  and on $\mb{x}_i$ (bounded fourth moment), covers the asymptotic normality of the LTS. 
V06c 
  also treated the general case $p\geq 1$ and obtained the asymptotic normality of the LTS under many stringent conditions on the non-random covariates $\mb{x}_i$s
and the distributions of $e_i$ in a twenty-seven pages long article. The assumption $\mathcal{C}$ is quite artificial and never verified there. 
On the other hand, C04 and C05 
also addressed the asymptotic normality of LTS in the nonlinear regression under a dependence setting. For these extensions, many artificial assumptions (D1, D2, H1, H2, H3, H4, H5, H6, I1, I2) are imposed, but they are never verified for the linear LTS case. So those results do not cover the  LTS in (\ref{lts.eqn}).
\vs
\tb{(c)}
Furthermore, since there was no population version like (\ref{objective.eqn}) and (\ref{lts-population.eqn}) before, empirical process theory could not be employed to verify the VC class of functions in V06a, b, c, and C04, C05.
 Our approach here is quite different from former classical analyses and much more neat and concise (the standard empirical process theory was asserted to be non-applicable in C04 and C05 due to the lack of a population version of the LTS).
\hfill \pend

\vs
\noindent
\section{Inference procedures} \label{Sec.4} 
In order to utilize the asymptotic normality result in Theorem 3.2, we need to figure out the asymptotic covariance. For simplicity, assume that $\bs{z}=(\bs{x}^{\top}, y)^{\top}$ follows  elliptical distributions $E(g; \bs{\mu}, \bs{\Sigma})$ with density
$$
f_{\bs{z}}(\bs{x}^{\top}, y)=\frac{g(((\bs{x}^{\top},y)^{\top}-\bs{\mu})^{\top}\bs{\Sigma}^{-1}((\bs{x}^{\top},y)^{\top}-\bs{\mu}))}{\sqrt{\det(\bs{\Sigma})}},
$$
where 
$\bs{\mu}\in \R^p$ and $\bs{\Sigma}$ is  a positive definite matrix of size $p$ which is proportional to the covariance matrix if the latter exists. We assume 
the $f_{\bs{z}}$ is unimodal.\vs

\vs
\noindent
\tb{Equivariance}
\vs
A regression estimation functional  $\bs{t}(\cdot)$ is said to be \emph{regression}, \emph{scale}, and \emph{affine} equivariant (see Zuo (2021a)) 
 if respectively
\bee
 \bs{t}(F_{(\bs{w},~ y+\bs{w}^{\top}b)})&=&\bs{t}(F_{(\bs{w},~y)})+\bs{b}, \forall ~\bs{b} \in \R^p;\\[1ex]
 \bs{t}(F_{(\bs{w},~sy)})&=&s \bs{t}(F_{(\bs{w},~y)}), \forall ~s \in \R;  \\[1ex]
 \bs{t}(F_{(\bs{A}^{\top}\bs{w}, ~y)})&=&\bs{A}^{-1}\bs{t}(F_{(\bs{w}, ~y)}),\forall \mbox{~nonsingular  $\bs{A}\in \R^{p\times p}$} ;
\ene
\vs
\noindent
\tb{Theorem 4.1} $\bs{\beta}_{lts}(F_{(\bs{x}', y)}, \alpha)$ is regression, scale, and affine equivariant.
\vs \noindent
\tb{Proof:}  See the empirical version treatment given in {RL87} (p. 132). \hfill \pend
\vs
\noindent
\tb{Transformation} ~
Assume the Cholesky decomposition of $\bs{\Sigma}$ yields a non-singular lower triangular matrix $\bs{L}$ of the form
\[
\left(
\begin{array}{cc}
\bs{A} & \bs{0}\\
\bs{v}^{\top}& c
\end{array}
\right)
\]
with $\bs{\Sigma}=\bs{L}\bs{L}^{\top}$. Hence $\det(\bs{A})\neq 0\neq c$. Now transfer $(\bs{x}^{\top}, y)$ to $(\bs{s}^{\top}, t)$ with  $(\bs{s}^{\top}, t)^{\top}=\bs{L}^{-1}((\bs{x}^{\top}, y)^{\top}-\bs{\mu})$. It is readily seen that the distribution of $(\bs{s}^{\top}, t)^{\top}$ follows
$E(g; \bs{0}, \bs{I_{p\times p}})$. \vs
Note that $(\bs{x}^{\top}, y)^{\top}=\bs{L}(\bs{s}^{\top}, t)^{\top}+(\bs{\mu}^{\top}_1, \mu_2)^{\top}$ with $\bs{\mu}=(\bs{\mu}^{\top}_1, \mu_2)^{\top}$. That is,
\begin{align}
\bs{x}&=\bs{A}\bs{s}+\bs{\mu}_1,\\[1ex]
y&=\bs{v}^{\top}\bs{s}+ct+\mu_2.
\end{align}
Equivalently,
\begin{align}
(1,\bs{s}^{\top})^{\top}&=\bs{B}^{-1}(1,\bs{x}^{\top})^{\top},   \label{x-transformation.eqn}\\[1ex]
t&=\frac{y-(1,\bs{s}^{\top})(\mu_2, \bs{v}^{\top})^{\top}}{c}, \label{y-transformation.eqn}
\end{align}
where
\[
\bs{B}=\begin{pmatrix}
1 &\bs{0}^{\top}\\
\bs{\mu}_1  &\bs{A}
\end{pmatrix}
,~~~~
\bs{B}^{-1}=
\begin{pmatrix}
1 &\bs{0}^{\top}\\
-\bs{A}^{-1}\bs{\mu}_1& \bs{A}^{-1}
\end{pmatrix},
\]
\vs
It is readily seen that (\ref{x-transformation.eqn}) is an affine transformation on $\bs{w}$ and (\ref{y-transformation.eqn}) is first an affine transformation on $\bs{w}$ then a regression transformation on $y$ followed by a scale transformation on $y$. In light of Theorem 4.1, we can assume hereafter, w.l.o.g. that $(\bs{x}^{\top}, y)$ follows an $E(g; \bs{0}, \bs{I}_{p\times p})$ (spherical) distribution and $\bs{I}_{p\times p}$ is the covariance matrix of $(\bs{x}^{\top}, y)$.
\vs
\noindent
\tb{Theorem 4.2} Assume that
 $e\sim\mathcal{N}(0, \sigma^2)$, $e$ and $\bs{x}$ are independent.
Then
\bi
\item[(1)] $P\nabla=\bs{0}$ and $P(\nabla\nabla')=8\sigma^2 C \bs{I}_{p\times p}$,\\[1ex]
 with $C=\Phi(c)-1/2-{c}e^{-c^2/2}/\sqrt{2\pi}$,
  where  $c=\sqrt{F^{-1}_{\chi^2(1)}(\alpha)}$ and $\Phi (x) $ is the CDF of $\mathcal{N}(0,1)$, $\chi^2(1)$ is a chi-square random variable with one degree of freedom.
\item[(2)] $\mb{V}= 2C_1\bs{I}_{p\times p} $ with $C_1=2\Phi(c)-1$.
\item[(3)] $n^{1/2}(\widehat{\bs{\beta}}^n_{lts} - \bs{\beta}_{lts}) \stackrel{d} \longrightarrow  {\cal{N}}(\bs{0}, \frac{2C\sigma^2}{C_1^2}\bs{I}_{p \times p})$,
where $C$ and $C1$ are defined in (1) and (2) above.
\ei
\vs \noindent
\tb{Proof:} See the Appendix. \hfill \pend
\vs
\noindent
\tb{Approximate $100(1-\gamma)\%$ confidence region}
\vs
\noin
\tb{(i) Based on the asymptotic normality}
~Under the setting of Theorem 4.2, an approximate $100(1-\gamma)\%$ confidence region for the unknown regression parameter $\bs{\beta}_0$ is:
$$
\Big\{\bs{\beta} \in \R^p:~~ \|\bs{\beta}-\widehat{\bs{\beta}}^n_{lts}\|\leq \sqrt{\frac{2C\sigma^2}{C_1^2n}}F^{-1}_{\chi^2(p)}(\gamma)\Big \},
$$
where $\|\cdot\|$ stands for the Euclidean distance. Without asymptotic normality, one can appeal to the next procedure.
\vs
\noin
\tb{(ii) Based on bootstrapping scheme and depth-median and depth-quantile} Here no assumptions on the underlying distribution are needed.
This approximate procedure first re-samples $n$ points with replacement from the given original sample points and  calculates an $\widehat{\bs{\beta}}^n_{lts}$. 
Repeat this $m$ (a large number, say $10^4$) times and obtain $m$ such  $\widehat{\bs{\beta}}^n_{lts}$s.\vs The next step is to calculate the depth, with respect to a location depth function (e.g., halfspace depth (Zuo (2018)) 
or projection depth (Zuo (2003) 
and Shao et al. (2022)), 
of these $m$ points in the parameter space of $\bs{\beta}$.
Trimming $\lfloor\gamma m \rfloor$ of the least deepest points among the $m$ points, the left points form a convex hull, that is an approximate  $100(1-\gamma)\%$ confidence region for the unknown regression parameter $\bs{\beta}_0$ (see Zuo (2009, 2010)) 
 in the location  and low dimensions cases).

\section{Concluding remarks} \label{Sec.5}

Without the population version for the LTS (see (\ref{lts-population.eqn})), it will be difficult to
apply the empirical process theory to study the asymptotics of the LTS, e.g, to verify the key
result, the VC-class property of regression function class (indexed by $\bs{\beta}$) will be challenging. 
Without an explicit regression function (unlike the linear case), one can
avoid this challenge, and simply
  assume this VC-class property.  This is what exactly done in nonlinear setting LTS in
 C04, C05. The latters even 
 believed that the standard empirical process theory is not applicable for the asymptotics of LTS 
 while V06a, b, c, addressed the asymptotics, without any advanced tools, employed elementary tools with numerous artificial assumptions 
 and lengthy articles.\vs
By partitioning the parameter space and introducing the population version for the LTS, this article establishes some fundamental and primary properties for the objective function of the LTS in both empirical and population settings. 
These newly obtained original results verify some key facts 
and facilitate the application of standard empirical process theory to 
 the establishment of asymptotic normality for the sample LTS in a concise and neat fashion. 
Some of newly obtained results, such as Fisher and strong consistency, and influence function,  are original and obtained as by-products. \vs
The asymptotic normality is applied in Theorem 4.2 for the practical inference procedure of confidence regions of the regression parameter $\bs{\beta}_0$. There are open problems left here; one is the estimation of the variance of $e$, which is now unrealistically assumed to be known, and the other is the testing of hypothesis on $\bs{\beta}_0$.
\vs
\vs
\noin
{\textbf{\Large Acknowledgments}}
\vs
Insightful comments and useful suggestions from Profs. Wei Shao and Derek Young  
have significantly improved the manuscript and are highly appreciated. Special thanks go to Prof. Derek Young for making the technical report of Chen, Stromberg, and Zhou available.

\vs
\noin
{\textbf{\Large Declarations}}
\vs
\noin
\tb{Funding}

This author declares that there is no funding received for this study.

\vs
\noin
\tb{Conflicts of interests/Competing interests}

This author declares that there is no conflict of interests/Competing interests.

\vs
\noindent
\tb{\Large Appendix: Proofs}
\vs
\noin
\tb{Proof of Lemma 2.1}

\vs
\noin
\tb{Proof}: 
(i) Based on the definition in (\ref{element-set.eqn}), over $S_{\bs{\beta}^l}$, there is no tie among the smallest $h$ squared residuals.  The assertion \tb{(a)} follows straightforwardly.
\vs
The first and the last equality in \tb{(b)} is trivial, it suffices to focus on the middle one.
Let $k_i:=k_i(\bs{\eta})$ ($=k_i(\bs{\beta}^l)$ in light of (\ref{element-set.eqn})).
By (\ref{objective-0.eqn}) we have that
\be O^n(\bs{\eta})=\frac{1}{n}\sum_{i=1}^h r^2_{k_i}(\bs{\eta}). \label{O-n-beta.eqn}
\ee
Let $r_i:=r_i(\bs{\eta})=y_i-\bs{w}^{\top}_i\bs{\eta}$ and
 $\gamma=\min\{\min_{1\leq i\neq j\leq n}\{|r^2_i-r^2_j|\}, 1\}$. 
 Then $1\geq \gamma>0$ (a.s.).
 \vs
 Based on the continuity of $r^2_{k_i}(\bs{\beta})$ in $\bs{\beta}$,  for any $1\leq i\leq h$ and for any given $\varepsilon\in (0, 1)$, we can fix a small $\delta>0$ so that $|r^2_{k_i}(\bs{\beta})-r^2_{k_i}(\bs{\eta})|<{\gamma\varepsilon}/{4h}$ for any $\bs{\beta}\in B(\bs{\eta},\delta)$.
   Now we have for any $\bs{\beta}\in B(\bs{\eta}, \delta)$, (assume below $2\leq i\leq h$ )
 \begin{align*}
 r^2_{k_i}(\bs{\beta})-r^2_{k_{(i-1)}}(\bs{\beta}) &> r^2_{k_i}(\bs{\eta})-\frac{\gamma\varepsilon}{4h}-[r^2_{k_{(i-1)}}(\bs{\eta})+\frac{\gamma \varepsilon}{4h}]\\[2ex]
 &= r^2_{k_i}(\bs{\eta})-r^2_{k_{(i-1)}}(\bs{\eta})-\frac{\gamma \varepsilon}{2h}\\[2ex]
 &\geq \gamma-\frac{\gamma \varepsilon}{2h}>0~ (a.s.),
\end{align*}
Thus, 
$k_i=k_i(\bs{\eta})$, $1\leq i\leq h$ forms the $h$-integer set
for any $\bs{\beta}\in B(\bs{\eta}, \delta)$.  Part \tb{(b)} follows.
\vs
(ii) The domain of function $O^n(\bs{\beta})$ is the union of the pieces of $\overline{S}_{\bs{\beta}^l}$ and the function of $O^n(\bs{\beta})$ over ${S}_{\bs{\beta}^l}$ is a quadratic function of $\bs{\beta}$:
$O^n(\bs{\beta})=\sum_{j=1}^hr^2_{k_j(\bs{\beta}^l)}(\bs{\beta})$, the statement follows.
\vs
(iii) By (ii), it is clear that $O^n(\bs{\beta})$ is continuous in $\bs{\beta}$ over each piece of $S_{\bs{\beta}^l}$. We only need to show that this holds true for any $\bs{\beta}$ that is on the boundary of $S_{\bs{\beta}^l}$.
\vs
Let $\bs{\eta}$ lie on the common boundary of $S_{\bs{\beta}^{s}}$ and $S_{\bs{\beta}^{t}}$, then $O^n(\bs{\beta})=\frac{1}{n}\sum_{i=1}^hr^2_{k_i(\bs{\beta}^s)}$ for any $\bs{\beta}\in \overline {S}_{\bs{\beta}^s}$
[this is obviously true if $\bs{\beta} \in S_{\bs{\beta}^{s}}$, it is also true if $\bs{\beta} $ on the boundary of $ S_{\bs{\beta}^{s}}$ since in this case the $\bs{\beta}$-$h$-integer set is not unique,  there are at least two, the one of them is $k_1(\bs{\beta}^s), \cdots, k_h(\bs{\beta}^s)$]
 and
$O^n(\bs{\beta})=\frac{1}{n}\sum_{i=1}^hr^2_{k_i(\bs{\beta}^t)}$ for any $\bs{\beta}\in \overline {S}_{\bs{\beta}^t}$. Let $\{\bs{\beta}_j\}$ be a sequence approaching  $\bs{\eta}$, where $\bs{\beta}_j$ could be on $\overline{S}_{\bs{\beta}^{s}}$ or on $\overline{S}_{\bs{\beta}^{t}}$. We show that $O^n(\bs{\beta}_j)$ approaches to $O^n(\bs{\eta})$. 
Note that $O^n(\bs{\eta})=\frac{1}{n}\sum_{i=1}^hr^2_{k_i(\bs{\beta}^s)}(\bs{\eta})=\frac{1}{n}\sum_{i=1}^hr^2_{k_i(\bs{\beta}^t)}(\bs{\eta})$. Partition $\{\bs{\beta}_j\}$ into $\{\bs{\beta}_{j_s}\}$ and $\{\bs{\beta}_{j_t}\}$ so that all members of the former belong to $\overline {S}_{\bs{\beta}^s}$ where the latter are all within $\overline{S}_{\bs{\beta}^{t}}$.
By continuity of the sum of $h$ squared residuals in $\bs{\beta}$, both $O^n(\bs{\beta}_{j_s}))$ and $O^n(\bs{\beta}_{j_t}))$ approach to $O^n(\bs{\eta})$ since both $\{\bs{\beta}_{j_s}\}$ and $\{\bs{\beta}_{j_t}\}$ approach $\eta$ as $\min \{s, t\} \to \infty$.
\vs
(iv) Note that for any $l$, $1\leq l\leq L$, over $S_{\bs{\beta}^l}$, one has a least squares problem with $n$ reduced to $h$, $O^n(\bs{\beta})$ is a quadratic function and hence is twice differentiable and strictly convex in light of the following
\begin{align*}
n\frac{\partial}{\partial\bs{\beta}} O^n(\bs{\beta}) &= -2\sum_{i=1}^nr_i \mathds{1}_i\bs{w}_i=-2\bs{X}^{\top}_n\bs{D}\bs{R}, \\[1ex]
n\frac{\partial^2}{\partial\bs{\beta}^2} O^n(\bs{\beta})&=2\bs{X}_n^{\top}\bs{D}\bs{X}_n=2 \bs{X^*}_n^{\top} \bs{X^*}_n=2\sum_{i=1}^h\bs{w}_{k_i(\bs{\beta}^l)}\bs{w}^{\top}_{k_i(\bs{\beta}^l)},
\end{align*}
where $\bs{R}=(r_1, r_2, \cdots, r_n)^{\top}$, $\bs{D}=\mbox{diag}(\mathds{1}_i)$, 
 and $\bs{X^*}_n=\bs{D}\bs{X}_n$. Strict convexity follows from the positive definite of the Hessian matrix: $\frac{2}{n}\bs{X^*}_n^{\top}\bs{X^*}_n$ (an invertible matrix due to \tb{(A1)}, see (iii) in the proof of Theorem 2.1). \hfill \pend
\vs

\vs
\noin
\tb{Proof of Theorem 2.1}
\vs
\noindent
\tb{Proof}:
{(i)}
Over each ${S}_{\bs{\beta}^l}$, an open set, $O^n(\bs{\beta})$ is twice differentiable and strictly convex in light of given condition,
 hence it has a unique minimizer (otherwise, one can show that by openness and strictly convexity there is a third point in ${S}_{\bs{\beta}^l}$ that attains a strictly smaller objective value than the two minimizers). Since there are only finitely many ${S}_{\bs{\beta}^l}$,  the assertion follows if we can prove that the minimum does not reach at a boundary point of some ${S}_{\bs{\beta}^l}$.
\vs
 Assume it is otherwise. That is, $O^n(\bs{\beta})$ reaches its global minimum at point $\bs{\beta}_1$ which is a boundary point of $S_{\bs{\beta}^l}$ for some $l$. Assume that over $S_{\bs{\beta}^l}$,
 $O^n(\bs{\beta})$ attains its local minimum value at the unique point $\bs{\beta}_2$. Then, $O^n(\bs{\beta}_1)\leq O^n(\bs{\beta}_2)$. If equality holds then
we have the desired result (since there are points besides $\bs{\beta}_2$ in ${S}_{\bs{\beta}^l}$ which also attain the minimum value as $\bs{\beta}_2$, a contradiction). Otherwise, there is a point $\bs{\beta}_3$ in the small neighborhood of $\bs{\beta}_1$ so that $O^n(\bs{\beta}_3)\leq O^n(\bs{\beta}_1)+(O^n(\bs{\beta}_2)-O^n(\bs{\beta}_1))/2< O^n(\bs{\beta}_2)$. A contradiction appears. 
\vs

{(ii)} It is seen from \tb{(i)} that $O^n(\bs{\beta})$ is twice continuously differentiable, hence its first derivative evaluated at the global minimum must be zero. By \tb{(i)}, we have the equation (\ref{estimation.eqn}). \vs

{(iii)} This part directly follows from \tb{(ii)} and the invertibility of $\bs{M}_n$. The latter follows from \tb{(A1)} that implies that the p columns of matrix $\bs{X}_n$ are linearly independent and also implies that any h sub-rows of $\bs{X}_n$ has a full rank. 
\hfill
\pend
\vs\vs
\noindent
\tb{Proof of Lemma 2.2}
\vs
\noin
\tb{Proof:}
Denote the integrand in (\ref{objective.eqn}) as $G(\bs{\beta}):=(y-\bs{w}^{\top}\bs{\beta})^2\mathds{1}\left((y-\bs{w}^{\top}\bs{\beta})^2\leq q(\bs{\beta},\alpha)\right)$ for a given point $(\bs{x}^{\top}, y)\in \R^p$.
Write $G(\bs{\beta}):=(y-\bs{w}^{\top}\bs{\beta})^2\big(1-\mathds{1}\left((y-\bs{w}^{\top}\bs{\beta})^2> q(\bs{\beta},\alpha)\right)\big)$.
\vs
\tb{(i)} By the strictly monotonicity of $F_W$ around $q(\bs{\beta},\alpha)$, we have the continuity of the $q(\bs{\beta},\alpha)$. Consequently, $G(\bs{\beta})$ is obvious continuous and so is  $O(\bs{\beta})$ in $\bs{\beta}\in \R^p$. 
\vs
\tb{(ii)} For arbitrary points $(\bs{x}^{\top}, y)$ and $\bs{\beta}$ in $\R^p$, there are three cases for the relationship between the squared residual and its quantile: (a) $(y-\bs{w}^{\top}\bs{\beta})^2>q(\bs{\beta},\alpha)$
(b) $(y-\bs{w}^{\top}\bs{\beta})^2<q(\bs{\beta},\alpha)$ and (c) $(y-\bs{w}^{\top}\bs{\beta})^2=q(\bs{\beta},\alpha)$. Case (c) happens with probability zero, we thus skip this case and treat (a) and (b) only.
By the continuity in $\bs{\beta}$, there is a small neighborhood of $\bs{\beta}$: $B(\bs{\beta}, \delta)$, centered at $\bs{\beta}$ with radius $\delta$ such that (a) (or (b)) holds for all $\bs{\beta} \in B(\bs{\beta}, \delta)$. This implies that 
$$
\frac{\partial}{\partial \bs{\beta}}\mathds{1}\left((y-\bs{w}^{\top}\bs{\beta})^2> q(\bs{\beta},\alpha)\right)=\bs{0}, ~(a.s.)
$$
and
\[
\frac{\partial}{\partial \bs{\beta}}G(\bs{\beta})=-2(y-\bs{w}^{\top}\bs{\beta})\bs{w}\mathds{1}\left((y-\bs{w}^{\top}\bs{\beta})^2\leq q(\bs{\beta},\alpha)\right), ~(a.s).
\]
Hence, we have that
\[
\frac{\partial^2}{\partial \bs{\beta}^2}G(\bs{\beta})=2\bs{w}\bs{w}^{\top}\mathds{1}\left((y-\bs{w}^{\top}\bs{\beta})^2\leq q(\bs{\beta},\alpha)\right), ~(a.s).
\]
Note that $E(\bs{w}\bs{w}^{\top})$ exists. 
 Then, by the Lebesgue dominated convergence theorem, the desired result follows.
\vs
\tb{(iii)} The strict convexity follows from the twice differentiability and the positive definite of the second order derivative of $O(\bs{\beta})$.
\hfill \pend

\vs
\vs
\noindent
\tb{Proof of Theorem 2.3}
\vs
\noin
\tb{Proof:}  We will treat $\bs{\beta}_{lts}(F_{(\bs{x}^{\top},~ y)},\alpha)$, the counterpart for $\bs{\beta}_{lts}(F_{\varepsilon}(\bs{z}),\alpha)$  can be treated analogously.
\vs
\tb{(i)} Existence follows from the positive semi-definiteness of the Hessian matrix (see proof of (ii) of Lemma 2.2) and the convexity of $O(\bs{\beta})$.\vs
\tb{(ii)} The equation follows from the differentiability and the first order derivative of $O(\bs{\beta})$ given in the proof (ii) of Lemma 2.2.\vs
\tb{(iii)} The uniqueness follows from the positive definite of the Hessian matrix based on the given condition (invertibility).
 \hfill \pend
 \vs \vs
\noindent
\tb{Proof of Theorem 2.4}
\vs
\noin
\tb{Proof:}  By theorem 2.3, (i) and given conditions guarantee the existence and the uniqueness of $\bs{\beta}_{lts}(F_{(\bs{x}^{\top}, y)}, \alpha)$ which is the unique solution of the system of the equations
\[
 \int(y-\bs{w}^{\top}\bs{\beta})\bs{w}\mathds{1}\left((y-\bs{w}^{\top}\bs{\beta})^2\leq q(\bs{\beta}, \alpha) \right) dF_{(\bs{x}^{\top}, y)}(\bs{x}, y)=\mb{0}.
 \]
 Notice that $y-\bs{w}^{\top}\bs{\beta}=-\bs{w}^{\top}(\bs{\beta}-\bs{\beta}_0)+e$. Inserting this into the above equation we have
 \[
 \int(-\bs{w}^{\top}(\bs{\beta}-\bs{\beta}_0)+e)\bs{w}\mathds{1}\left((-\bs{w}^{\top}(\bs{\beta}-\bs{\beta}_0)+e)^2\leq F^{-1}_{(-\bs{w}^{\top}(\bs{\beta}-\bs{\beta}_0)+e)^2}(\alpha)\right) dF_{(\bs{x}^{\top}, y)}(\bs{x}, y)=\mb{0}.
 \]
 By (ii) it is readily seen that $\bs{\beta}=\bs{\beta}_0$ is a solution of the above system of equations. Uniqueness leads to the desired result.
 \hfill \pend
 \vs \vs
\noindent
\tb{Proof of Theorem 2.5}
\vs
\noin
\tb{Proof:} Write $\bs{\beta}^{\varepsilon}_{lts}(\mb{z}_0)$ for $\bs{\beta}_{lts}(F_{\varepsilon}(\mb{z}_0), \alpha)$ and insert it for $\bs{\beta}$  into (\ref{lts-contamination-estimation.eqn}) and take derivative with respect to $\varepsilon$ in both sides of (\ref{lts-contamination-estimation.eqn}) and let $\varepsilon \to 0$, we obtain (in light of dominated theorem)
 \begin{align}
& \left(\int \frac{\partial}{\partial \bs{\beta}^{\varepsilon}_{lts}(\mb{z}_0)} r(\bs{\beta}^{\varepsilon}_{lts}(\mb{z}_0))\mb{v}\mathds{1}\left( r(\bs{\beta}^{\varepsilon}_{lts}(\mb{z}_0))^2\leq q_{\varepsilon}(\mb{z},\bs{\beta}^{\varepsilon}_{lts}(\mb{z}_0),\alpha) \right) \Bigg|_{\varepsilon \to 0}\!\!\! dF_{(\bs{x}^{\top} y)}\right) \mbox{IF}(\mb{z}_0, \bs{\beta}_{lts}, F_{(\bs{x}^{\top},y )})\nonumber \\[2ex]
&+ \int r(\bs{\beta}_{lts}(F_{(\mb{x}^{\top}, y)}, \alpha)))\bs{w}\mathds{1}\left(r(\bs{\beta}_{lts}(F_{(\mb{x}^{\top}, y)}, \alpha))))^2\leq q(\bs{\beta}_{lts}(F_{(\mb{x}^{\top}, y)}, \alpha), \alpha) \right) d(\delta_{\mb{z}_0}-F_{(\mb{x}^{\top}, y)}) \nonumber \\[2ex]
&=\mb{0}, \label{if-proof.eqn}
 \end{align}
 where $r(\bs{\beta})=t-\bs{v}^{\top}\bs{\beta}$ in the first term on the LHS and $r(\bs{\beta})=y-\bs{w}^{\top}\bs{\beta}$ in the second term on the LHS.
 Call the two terms on the LHS as  $T_1$ and $T_2$ respectively and call the integrand in $T_1$ as $T_0$, then it is seen that (see the proof \tb{(i)} of Theorem 2.1)
 \bee
 T_0&=&\frac{\partial}{\partial \bs{\beta}^{\varepsilon}_{lts}(\mb{z}_0)} (t-\bs{v}^{\top}\bs{\beta}^{\varepsilon}_{lts}(\mb{z}_0))\mb{v}\mathds{1}\left((t-\bs{v}^{\top}\bs{\beta}^{\varepsilon}_{lts}(\mb{z}_0))^2\leq q_{\varepsilon}(\mb{z},\bs{\beta}^{\varepsilon}_{lts}(\mb{z}_0),\alpha) \right) \Bigg|_{\varepsilon \to 0}\\[1ex]
 &=&-\bs{w}\bs{w}^{\top}\mathds{1}\left((y-\bs{w}^{\top}\bs{\beta}_{lts})^2\leq q(\bs{\beta}_{lts},\alpha) \right).
 \ene
 Focus on the $T_2$, it is readily seen that
 \begin{align*}
 T_2&=\int (y-\bs{w}^{\top}\bs{\beta}_{lts}(F_{(\mb{x}^{\top}, y)}, \alpha))\bs{w}\mathds{1}\left((y-\bs{w}^{\top}\bs{\beta}_{lts}(F_{(\mb{x}^{\top}, y)}, \alpha))^2\leq q(\bs{\beta}_{lts}(F_{(\mb{x}^{\top}, y)}, \alpha), \alpha) \right) d\delta_{\mb{z}_0}\\[1ex]
 &-\int (y-\bs{w}^{\top}\bs{\beta}_{lts}(F_{(\mb{x}^{\top}, y)}, \alpha))\bs{w}\mathds{1}\left((y-\bs{w}^{\top}\bs{\beta}_{lts}(F_{(\mb{x}^{\top}, y)}, \alpha))^2\leq q(\bs{\beta}_{lts}(F_{(\mb{x}^{\top}, y)}, \alpha), \alpha) \right) dF_{(\mb{x}^{\top}, y)},
 \end{align*}
 In light of (\ref{lst-estimation.eqn}) we have
 \begin{align*}
 T_2&=\int (y-\bs{w}^{\top}\bs{\beta}_{lts}(F_{(\mb{x}^{\top}, y)}, \alpha))\bs{w}\mathds{1}\left((y-\bs{w}^{\top}\bs{\beta}_{lts}(F_{(\mb{x}^{\top}, y)}, \alpha))^2\leq q(\bs{\beta}_{lts}(F_{(\mb{x}^{\top}, y)}, \alpha), \alpha) \right) d\delta_{\mb{z}_0}\\
 &=\left\{
 \begin{array}{ll}
 \mb{0}, & \mbox{if $(t_0-\bs{v}^{\top}_0\bs{\beta}_{lts})^2>q(\bs{\beta}_{lts}, \alpha)$,}\\[2ex]
 (t_0-\bs{v}^{\top}_0 \bs{\beta}_{lts})\bs{v}_0,&    \mbox{otherwise}.
 \end{array}
 \right.
 \end{align*}
 This,  $T_0$, and display (\ref{if-proof.eqn}) lead to the desired result.
 \hfill \pend

 \vs\vs
\noindent
\tb{Proof of Lemma 3.1}
\vs
\noin
\tb{Proof}: We invoke Theorem 24 in II.5 of Pollard (1984) (P84). 
The first requirement of the theorem is the existence of an envelope of  $ \mathscr{F}$. The latter is $\sup_{\bs{\beta} \in \Theta} F^{-1}_{r(\bs{\beta})^2}(c)$, which is bounded since $\Theta$ is compact and $F^{-1}_W$ is continuous in $\bs{\beta}$, and $F^{-1}_W(\alpha)$ is non-decreasing in $\alpha \in [1/2, c]$. 
To complete the proof, we only need to verify the second requirement of the theorem.\vs
For the second requirement, that is, to bound the covering numbers, it suffices to show that
the graphs of functions in $\mathscr{F}(\bs{\beta}, \alpha)$ have only polynomial discrimination (see Theorem 25 and Example 26 in II.5 of {P84}).
\vs
The graph of a real-valued function $f$ on a set $S$ is defined as the subset (see p. 27 of {P84})
$$G_f = \{(s, t): 0\leq t \leq f(s) ~\mbox{or}~ f(s)\leq t \leq 0, s \in S \}.$$
\vs
The graph of a function in $\mathscr{F}(\bs{\beta}, \alpha)$ contains a point $(\mb{x^{\top}(\omega)}, y(\omega), t)$  if and only if
$0\leq t\leq f(\bs{x}, y, \bs{\beta}, \alpha)$ or $ f(\bs{x},y, \bs{\beta}, \alpha) \leq t\leq 0$. The latter case could be excluded since the function is always nonnegative (and equals $0$ case covered by the former case). The former case happens if and only if
$0\leq\sqrt{t}\leq y-\bs{w^{\top}}\bs{\beta}$ or $0\leq\sqrt{t}\leq -y+\bs{w^{\top}}\bs{\beta}$.
\vs
Given a collection of $n$ points $(\bs{x}^{\top}_i, y_i, t_i)$ ($t_i\geq 0$),  the graph of a function in $\mathscr{F}(\bs{\beta}, \alpha)$ picks out only points that belong  
to $\{\sqrt{t_i}\geq 0\}\cap \{y_i-\bs{\beta}^{\top}\bs{w}_i-\sqrt{t_i}\geq 0\} \cup \{\sqrt{t_i}\geq 0\}\cap\{-y_i+\bs{\beta}^{\top}\bs{w}_i-\sqrt{t_i}\geq 0\} $. 
Introduce $n$ new points $(\bs{w}^{\top}_i, y_i, z_i):=((1,\bs{x}^{\top}_i), y_i, \sqrt{t_i})$ in $\R^{p+2}$. On $\R^{p+2}$ define a vector space $\mathscr{G}$ of functions
$$g_{a, b, c}(\bs{w}, y, z)=\mb{a}^{\top}\bs{w}+by+cz,$$
where $a\in \R^p$, $b\in \R^1$, and $c\in \R^1$ and $\mathscr{G}:=\{g_{a, b, c}(\bs{w}, y, z)=\mb{a}^{\top}\bs{w}+by+cz, a\in \R^p, b\in \R^1, ~\mbox{and}~ c\in \R^1\}$ which is a $\R^{p+2}$-dimensional vector space.
\vs

It is clear now that the graph of a function in $\mathscr{F}(\bs{\beta}, \alpha)$ picks out only points that belong to
 the sets of $\{g\geq 0\}$ for $g\in \mathscr{G}$ (ignoring the union and intersection operations at this moment). By Lemma 18 in II.4 of {P84} (p. 20), the graphs of functions in $\mathscr{F}(\bs{\beta}, \alpha)$  pick only polynomial numbers of subsets of $\{\bs{p}_i:=(\bs{w}^{\top}_i, y_i, z_i), i\in\{1,\cdots, n\}\}$; those sets corresponding to
$g\in  \mathscr{G}$ with $a \in \{\bs{0},-\bs{\beta}, \bs{\beta}\}$, $b\in \{0, 1, -1\}$, and $c \in \{1, -1\}$  pick up even few subsets from $\{\bs{p}_i, i\in\{1,\cdots, n\}\}$. This in conjunction with Lemma 15 in II.4 of {P84} (p. 18), yields that
the graphs of functions in $\mathscr{F}(\bs{\beta}, \alpha)$  have only polynomial discrimination. \vs
By Theorem 24 in II.5 of {P84} we have completed the proof.
\hfill \pend
 \vs \vs
\noindent
\tb{Proof of Theorem 3.2}
\vs
\noin
\tb{Proof}: To apply Lemma 3.2, we need to verify the five conditions, among them only (iii) and (v) need to be addressed, all others are satisfied trivially. For (iii), it holds automatically since our $\tau_n=\widehat{\bs{\beta}}^n_{lts}$ is defined to be the minimizer of  $F_n(t)$ over $t\in T (=\Theta)$.\vs
So the only condition that needs to be verified is the (v), the stochastic equicontinuity of $\{E_nr(\cdot, t)\}$ at $t_0$. For that, we will appeal to the Equicontinuity Lemma (VII.4 of {P84}, p. 150).
To apply the Lemma, we will verify the condition for the  random covering numbers satisfy the uniformity condition. To that end, we look at the class of functions
\[
\mathscr{R}(\bs{\beta}, \alpha)=\left\{r(\cdot, \cdot, \alpha, \bs{\beta})=\left(\frac{\bs{\beta}^{\top}}{\|\bs{\beta}\|}V/2 \frac{\bs{\beta}}{\|\bs{\beta}\|} \right)\|\bs{\beta}\|:~ \bs{\beta}\in \Theta, \alpha\in [1/2,c]\right\}.
\]
Obviously, $\lambda_{max} r_0/2$ is an envelope for the class $\mathscr{R}$ in $\mathscr{L}^2(P)$, where $r_0$ is the radius of the ball $\Theta=B(\bs{\beta}_{lts}, r_0)$. We now show that the covering numbers of $\mathscr{R}$ is uniformly bounded, which amply suffices for the Equicontinuity Lemma.  For this, we will invoke Lemmas II.25 and II.36 of {P84}. 
 To apply Lemma II.25, we need to show that the graphs of functions in $\mathscr{R}$ have only polynomial discrimination. The graph of $r(\bs{x}, y, \alpha, \bs{\beta})$ contains a point $(\bs{x}^{\top}, y, t)$, $t\geq 0$
  if and only if $\left(\frac{\bs{\beta}^{\top}}{\|\bs{\beta}\|}V/2 \frac{\bs{\beta}}{\|\bs{\beta}\|} \right)\|\bs{\beta}\| \geq t$ for all $\bs{\beta} \in \Theta$ and $\alpha \in [1/2, c]$.
  \vs
 Equivalently, the graph of $r(\bs{x}, y, \alpha, \bs{\beta})$ contains a point $(\bs{x}^{\top}, y, t)$, $t\geq 0$ if and only if $\lambda_{min}/2\|\bs{\beta}\|\geq t$.
For a collection of $n$ points $(\bs{x}^{\top}_i, y_i, t_i)$ with $t_i\geq 0$, the graph picks out those points satisfying $\lambda_{min}/2\|\bs{\beta}\|- t_i\geq 0$. Construct from
$(\bs{x}^{\top}_i, y_i, t_i)$ a point $z_i=t_i$ in $\R$. On $\R$ define a vector space $\mathscr{G}$ of functions
\[
g_{a, b}(x)=ax+b,~~ a,~ b \in \R.
\]
By Lemma 18 of {P84}, the sets $\{g\geq 0\}$, for $g \in \mathscr{G}$, pick out only a polynomial number of subsets from $\{z_i\}$; those sets corresponding to functions in $\mathscr{G}$ with $a=-1$ and $b=\lambda_{min}/2\|\bs{\beta}\|$ pick out even fewer subsets from  $\{z_i\}$. Thus the graphs of functions in $\mathscr{R}$ have only    polynomial discrimination.
\hfill \pend

\vs\vs
\noindent
\tb{Proof of Theorem 4.2}: \vs
\noin
\tb{Proof}: In order to invoke Theorem 3.2, we only need to check the uniqueness of $\widehat{\bs{\beta}}^n_{lts}$ and ${\bs{\beta}}_{lts}$. The former is guaranteed by the (iii) of Theorem 2.1 since \tb{(A1)} holds true a.s.. This is because that any p columns of the
$\bs{X}_n$ or any its h rows could be regarded as a sample from an absolutely continuous random vector with dimension $n$ or $h$. The probability that these p points lie in a $(p-1)$ dimensional non-degenerated hyperplane (the normal vector is non-zero) is zero. \vs The latter is guaranteed by the (iii) of Theorem 2.3 since $W=(y-\bs{w}^{\top}\bs{\beta})^2$ is the square of a normal distribution with mean $-\beta_1$ and hence has a positive density and furthermore (18) becomes $2(\Phi(c/\sigma)-1/2)\bs{I}_{p\times p}$ hence is invertible, where $c$ is defined in Theorem 4.2.
By Theorems 4.1 and 2.4, we can assume, w.l.o.g., that ${\bs{\beta}}_{lts}=\bs{\beta}_0=\bs{0}$. Utilizing the independence between $e$ and $\bs{x}$ and Theorem 3.2, a straightforward calculation leads to the results.
\hfill \pend
\end{document}